\documentclass[12pt,a4paper]{article}
\usepackage{amsmath,amsfonts,amsthm,amsopn,color,amssymb}
\usepackage{graphicx}
\newcommand{\eps}{\varepsilon}
\newcommand{\R}{\mathbb{R}}

\newcommand{\RN}{{\mathbb{R}^N}}
\newcommand{\RP}{{\mathbb{R}^N_+}}
\newcommand{\de}{\partial}

\DeclareMathOperator{\cat}{cat}
\DeclareMathOperator{\dv}{div}

\renewcommand{\le}{\leqslant}
\renewcommand{\ge}{\geqslant}
\renewcommand{\a }{\alpha }

\renewcommand{\b }{\beta }

\renewcommand{\d }{\delta }

\newcommand{\g }{\gamma }

\renewcommand{\l }{\lambda}
\newcommand{\n }{\nabla }

\renewcommand{\O}{\Omega}
\renewcommand{\OE}{\Omega_\varepsilon}
\newcommand{\G}{\Gamma}

\renewcommand{\S}{\Sigma}
\renewcommand{\L}{\Lambda}

\newcommand{\A}{{\cal A}}

\renewcommand{\P}{{\cal P}}

\newcommand{\N}{\mathbb{N}}

\newcommand{\tf}{\tilde{\phi_1}}

\newtheorem{theorem}{Theorem}[section]
\newtheorem{lemma}[theorem]{Lemma}
\newtheorem{example}[theorem]{Example}
\newtheorem{definition}[theorem]{Definition}
\newtheorem{proposition}[theorem]{Proposition}
\newtheorem{remark}[theorem]{Remark}

\renewenvironment{proof}{\noindent{\textbf{Proof\quad}}}{$\hfill\square$\vspace{0.2 cm}\\}
\newenvironment{proof2}{\noindent{\textbf{Proof of Theorem \ref{th2}\quad}}}{$\hfill\square$\vspace{0.2 cm}\\}

\begin{document}

\title{\textbf{Singularly perturbed Neumann problems with potentials}}
\author{Alessio Pomponio\thanks{Supported by MIUR, national project \textit{Variational methods and nonlinear differential equations} } 
\\ 
\
\\
SISSA, via Beirut 2/4 \\ I-34014 Trieste 
\\ 
{\it pomponio@sissa.it}}

\date{ }

\maketitle

\section{Introduction}

In this paper we study the following problem:
\begin{equation}\label{EQe}
\left\{
\begin{array}[c]{ll}
-\eps^2 \dv \left(J(x)\nabla u\right)+V(x)u=u^p 
& \text{in }\O,
\\
\dfrac{\de u}{\de \nu}=0
& \text{on }\de\O,
\end{array}
\right.  
\end{equation}
where $\O$ is a smooth bounded domain with external normal $\nu$, $N\ge 3$, $1<p<(N+2)/(N-2)$, 
$J\colon\RN \to \R$ and $V\colon\RN \to \R$ are $C^2$ functions.

When $J\equiv 1$ and $V\equiv 1$, then \eqref{EQe} becomes 
\begin{equation}\label{EQN}
\left\{
\begin{array}[c]{ll}
-\eps^2 \varDelta u +u=u^p 
& \text{in }\O,
\\
\dfrac{\de u}{\de \nu}=0
& \text{on }\partial\O.
\end{array}
\right.  
\end{equation}
Such a problem was intensively studied in several works. For example, Ni \& Takagi, in \cite{NT, NT2}, 
show that, for $\eps$ sufficiently small, there exists a solution $u_\eps$ of \eqref{EQN} 
which concentrates in a point $Q_\eps \in \de \O$ and  $H(Q_\eps) \to \max_{\de \O} H$, here $H$ 
denotes the mean curvature of $\de \O$. Moreover in \cite{Li}, using the Liapunov-Schmidt reduction, 
Li constructs solutions with single peak and multi-peaks on $\de \O$ located near any stable 
critical points of $H$. Since the publication of \cite{NT, NT2}, there have been many works 
on spike-layer solutions of \eqref{EQN}, see for example \cite{DY, dPFW, GPW, G, GW, W2} and 
references therein.

What happens in presence of potentials $J$ and $V$?

In this paper we try to give an answer to this question and we will show that, 
for the existence of concentrating solutions, one has to check if at least one 
between $J$ and $V$ is not constant on $\de \O$. In this case the concentration 
point is determined by $J$ and $V$ only. In the other case the concentration 
point is determined by an interplay among the derivatives of $J$ and $V$ calculated on $\de \O$ 
and the mean curvature $H$.

On $J$ and $V$ we will do the following assumptions: 
\begin{description}
\item[(J)] $J\in C^2 (\O, \R)$, $J$ and $D^2 J$ are bounded; moreover, 
\begin{equation*}
J(x)\ge C>0 \quad \textrm{for all } x\in\O;
\end{equation*}
\item[(V)] $V\in C^2 (\O, \R)$, $V$ and $D^2 V$ are bounded; moreover, 
\begin{equation*}
V(x)\ge C>0 \quad \textrm{for all } x\in\O.
\end{equation*}
\end{description}

Let us introduce an auxiliary function which will play a crucial r\^ole in the study of \eqref{EQe}. 
Let $\G\colon \de \O \to \R$ be a function so defined:
\begin{equation}\label{eq:Gamma}
\G(Q)=V(Q)^{\frac{p+1}{p-1}- \frac{N}{2}} J(Q)^{\frac{N}{2}}.
\end{equation}
Let us observe that by {\bf (J)} and {\bf (V)}, $\G$ is well defined.

Our first result is:

\begin{theorem}\label{th1}
Let  $Q_0 \in \de \O$. Suppose {\bf (J)} and {\bf (V)}. There exists $\eps_0>0$ such that 
if $0<\eps<\eps_0$, 
then \eqref{EQe} possesses a solution $u_\eps$ 
which concentrates in $Q_\eps$ with $Q_\eps \to Q_0$, as $\eps \to 0$, 
provided that one of the two following conditions holds:
\begin{description}
\item[$(a)$] $Q_0$ is a non-degenerate critical point of $\G$;
\item[$(b)$] $Q_0$ is an isolated local strict minimum or maximum of $\G$.
\end{description}
\end{theorem}

Hence, if $J$ and $V$ are not constant on the boundary $\de \O$, the concentration phenomena 
depend only by $J$ and $V$ and not by the mean curvature $H$. 
Our second result deals with the other case and, more precisely, we will show 
that, if $J$ and $V$ (and so also $\G$) are constant on the boundary, then 
the concentration phenomena are due by another auxiliary function 
which depends on the derivatives of $J$ and 
$V$ on the boundary and by the mean curvature $H$. Let $\bar \S\colon \de \O \to \R$ be the function so 
defined:
\begin{multline}\label{eq:Sigmabar}
\bar \S(Q)\equiv\
k_1 \int_{\R^-_{\nu(Q)}} J'(Q)[x] \left|\left(\n \bar U \right)\!\left(k_2 x \right)\right|^2 d x
\\
+k_3 \int_{\R^-_{\nu(Q)}} V'(Q)[x] \left[ \bar U \!\!\left(k_2 x\right)\right]^2 d x
-k_4 H(Q),
\end{multline}
where $\bar U$ is the unique solution of
\[
\left\{
\begin{array}
[c]{lll}
-\varDelta  \bar U+ \bar U= \bar U^p &   \text{in }\R^{N},
\\
\bar U>0 &   \text{in }\R^{N},
\\
\bar U(0)=\max_{\R^{N}}\bar U,
\end{array}
\right.
\]
$\nu(Q)$ is the outer normal in $Q$ at $\O$,
\[
\R^-_{\nu(Q)} \equiv \left\{ x\in\RN : x \cdot \nu(Q)  \le 0\right\},
\]
and, for $i=1, \ldots, 4$, $k_i$ are constants which depend only on $J$ and $V$ and not on $Q$ 
(see Remark \ref{re:Sigma} for an explicit formula).

Our second result is:

\begin{theorem}\label{th2}
Suppose {\bf (J)} and {\bf (V)} with $J$ and $V$ constant on the boundary $\de \O$. 
Let  $Q_0 \in \de \O$ be an isolated local strict minimum or maximum of $\bar \S$. 
There exists $\eps_0>0$ such that if $0<\eps<\eps_0$, 
then \eqref{EQe} possesses a solution $u_\eps$ 
which concentrates in $Q_\eps$ with $Q_\eps \to Q_0$, as $\eps \to 0$. 
\end{theorem}

\begin{example}
Suppose that $J\equiv 1$ and fix any $Q_0 \in \de \O$. For $k\in \N$, let $V_k$ be a 
bounded smooth function constantly 
equal to $1$ on the $\de \O$ and in the whole $\O$, except a little ball tangent 
at $\de \O$ in $Q_0$, with $\n V_k (Q_0)=-k \nu(Q_0)$  (see figure 1).

\begin{figure}
\begin{center}
\includegraphics[height=4.8cm]{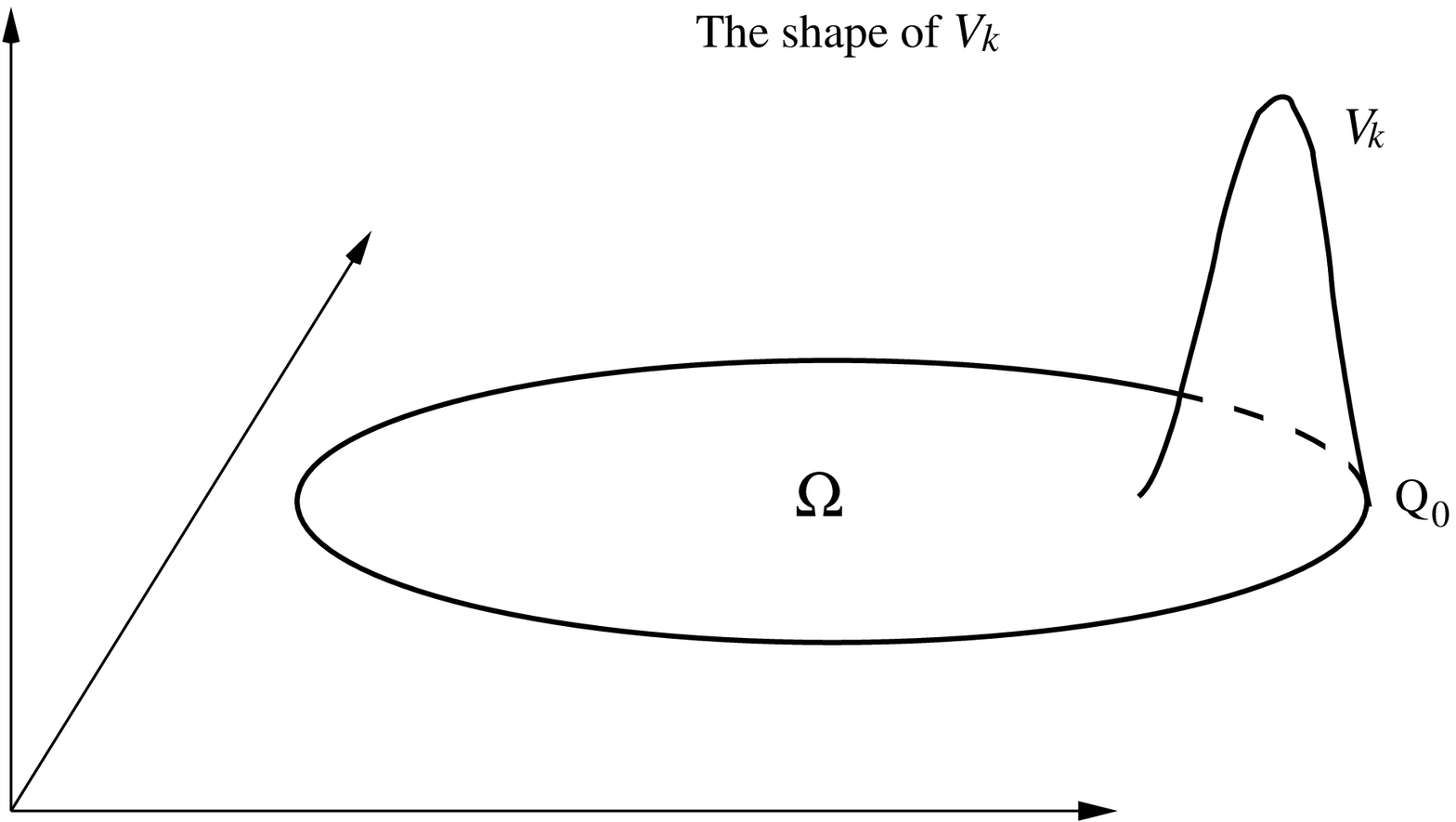}
\\ \hspace{0.2cm}\footnotesize{Figure 1}
\end{center}
\end{figure}

\noindent It is easy to see that, outside a little neighborhood of $Q_0$ in $\de \O$, we have 
\[
\bar \S(Q)= -C_1 H(Q),
\]
while
\[
\bar \S(Q_0)= -C_1 H(Q_0) + k C_2,
\]
where
\begin{eqnarray*}
C_1 &=& \frac{1}{2} \bar{B} + \left(\frac{1}{2}-\frac{1}{p+1}\right) \bar{A}, 
\\
C_2 &=& - \frac{1}{2}
\int_{\{ \nu(Q_0) \cdot x\le 0 \}} \nu(Q_0) \cdot x \;\bar U^2 d x.
\end{eqnarray*}
Since $C_2>0$, we can choose $k\gg 1$ such that $Q_0$ is the absolute maximum point for $\bar \S$ 
and hence there exists a solution concentrating at $Q_0$. 
\end{example}

Theorem \ref{th1} will be proved as a particular case of two multiplicity results in Section 6, where 
we will prove also Theorem \ref{th2}. The proof of the theorems relies on a finite dimensional 
reduction, precisely on the perturbation technique developed in \cite{AB, ABC, AMS}. 
In Section 2 we give some preliminary 
lemmas and some estimates which will be useful in Section 3 and Section 4, 
where we perform the Liapunov-Schmidt reduction, 
and in Section 5, where we make the asymptotic expansion of the finite dimensional functional.

Finally we mention that problem \eqref{EQe}, but with the Dirichlet boundary 
conditions, is studied by the author 
and by S. Secchi in \cite{PS}, where we  
show that there are solutions which concentrate in minima of an auxiliary 
function, which depends only on $J$ and $V$.

\

\noindent {\bf Acknowledgments} \quad The author wishes to thank Professor Antonio Ambrosetti and 
Professor Andrea Malchiodi for suggesting the problem and for useful discussions.

\begin{center}{\bf Notation}\end{center}
\begin{itemize}
\item $\R^N_+\equiv \left\{ (x_1, \ldots, x_N)\in\RN : x_N > 0\right\}$.

\item If $\mu \in \RN$, then 
$\R^-_\mu \equiv \left\{ x\in\RN :  x \cdot \mu  \le 0\right\}$, 
where with $x \cdot \mu$ we denote the scalar product in $\RN$ between $x$ and $\mu$.

\item If $r>0$ and $x_0 \in \RN$, $B_r (x_0)\equiv \left\{ x\in\RN : |x- x_0| <r \right\}$. 
We denote with $B_r$ the ball of radius $r$ centered in the origin.

\item If $u \colon \RN \to \R$ and $P\in \RN$, we set $u_P \equiv u(\cdot - P)$.

\item If $U^Q$ is the function defined in \eqref{eq:UQ}, when there is no misunderstanding, 
we will often write $U$ instead of $U^Q$. Moreover if $P=Q/\eps$, then $U_P \equiv U^Q(\cdot -P)$.

\item If $Q\in \de \O$, we denote with $\nu(Q)$ the outer normal in $Q$ at $\O$ and with $H(Q)$ the mean 
curvature of $\de \O$ in $Q$. 

\item If $\eps>0$, we set $\OE \equiv \O/\eps \equiv\{x\in \RN : \eps x \in \O\}$.

\item We denote with $\|\cdot \|$ and with $(\cdot \mid \cdot)$ respectively the norm and the scalar 
product of $H^1(\OE)$. While we denote with $\|\cdot \|_+$ and with $(\cdot \mid \cdot)_+$ 
respectively the norm and the scalar product of $H^1(\RP)$.

\item If $P \in \de \O_\eps$, we set $\de_{P_i} \equiv \frac{\de}{\de e_i}$, 
where $\left\{  e_1,\ldots,e_{N-1}\right\}$ is an orthonormal basis of $T_P(\partial\O_\eps)$. Analogously, 
if $Q\in \de \O$, we set $\de_{Q_i} \equiv \frac{\de}{\de {\tilde e}_i}$, 
where $\left\{  {\tilde e}_1,\ldots,{\tilde e}_{N-1}\right\}$ is an orthonormal basis of $T_Q(\de\O)$.
\end{itemize}

\section{Preliminary lemmas and some estimates}

First of all we perform the change of variable $x \mapsto \eps x$ and so problem \eqref{EQe} becomes
\begin{equation}\label{EQ}
\left\{
\begin{array}[c]{ll}
-\dv \left(J(\eps x)\nabla u\right)+V(\eps x)u=u^p 
& \text{in }\O_\eps,
\\
\dfrac{\de u}{\de \nu}=0
& \text{on }\de\O_\eps,
\end{array}
\right.  
\end{equation}
where $\O_\eps=\eps^{-1}\O$. Of course if $u$ is a solution of \eqref{EQ}, then $u(\cdot/\eps)$ 
is a solution of \eqref{EQe}.

Solutions of (\ref{EQ}) are critical points $u\in H^1(\O_\eps)$ of
\[
f_\eps (u)=
\frac{1}{2}\int_{\OE} J(\eps x)|\nabla u|^2 dx+
\frac{1}{2}\int_{\OE} V(\eps x)u^2dx
-\frac{1}{p+1}\int_{\OE} |u|^{p+1}.
\]
The solutions of (\ref{EQ}) will be found near a $U^Q$, the unique solution of
\[
\left\{
\begin{array}
[c]{lll}
-J(Q)\varDelta  u+ V(Q) u= u^p &   \text{in }\R^{N},\\
u>0 &   \text{in }\R^{N},\\
u(0)=\max_{\R^{N}} u, & 
\end{array}
\right.
\]
for an appropriate choice of $Q \in \de\O$. It is easy to see that 
\begin{equation}\label{eq:UQ}
U^Q(x)=V(Q)^{\frac{1}{p-1}}\,\bar U\left(x \sqrt{V(Q)/J(Q)}  \right),
\end{equation}
where $\bar U$ is the unique solution of
\[
\left\{
\begin{array}
[c]{lll}
-\varDelta  \bar U+ \bar U= \bar U^p &   \text{in }\R^{N},
\\
\bar U>0 &   \text{in }\R^{N},
\\
\bar U(0)=\max_{\R^{N}}\bar U,
\end{array}
\right.
\]
which is radially symmetric and decays exponentially at infinity with its derivatives.

We remark that  $U^Q$ is a solution also of the ``problem to infinity'':
\begin{equation}\label{eq:Q}
\left\{
\begin{array}
[c]{lll}
-J(Q) \varDelta  u+ V(Q) u= u^p &   \text{in }\R^{N}_+,
\\
\dfrac{\de u}{\de \nu}=0 & \text{on } \de\R^N_+.
\end{array}
\right.
\end{equation}
The solutions of (\ref{eq:Q}) are critical points of the functional defined on $H^1(\RP)$
\begin{equation}\label{eq:F}
F^{Q}(u)=
\frac{1}{2}J(Q)\int_{\R^N_+}|\nabla u|^2 +
\frac{1}{2}V(Q)\int_{\R^N_+}u^2
-\frac{1}{p+1}\int_{\RP}|u|^{p+1}.
\end{equation}

We recall that we will often write $U$ instead of $U^Q$. 
If $P=\eps^{-1}Q \in \de\OE$, we set $U_P \equiv U^Q(\cdot - P)$ and 
\[
Z^\eps \equiv \{ U_P : P\in \de\OE\}.
\]

\begin{lemma}\label{lem:nf}
For all $Q \in \de \O$ and for all $\eps$ sufficiently small, if $P=Q/\eps \in \de \OE$, then
\begin{equation}\label{eq:nf}
\|\nabla f_\eps(U_P)\|=O(\eps).
\end{equation}
\end{lemma}

\begin{proof}
\begin{gather*}
(\nabla f_\eps(U_P) \mid v) =
\int_{\OE} J(\eps x) \nabla U_P \cdot \n v
+\int_{\OE} V(\eps x) U_P v 
-\int_{\OE} U_P^p v
\\
=\int_{\frac{\O -Q}{\eps}} J(\eps x +Q) \nabla U \cdot \n v_{-P}
+\int_{\frac{\O -Q}{\eps}} V(\eps x +Q) U v_{-P} 
-\int_{\frac{\O -Q}{\eps}} U^p v_{-P}
\\
=\int_{\frac{\O -Q}{\eps}} J(Q) \nabla U \cdot \n v_{-P}
+\int_{\frac{\O -Q}{\eps}} V(Q) U v_{-P} 
-\int_{\frac{\O -Q}{\eps}} U^p v_{-P}
\\
+\int_{\frac{\O -Q}{\eps}} (J(\eps x +Q)-J(Q)) \nabla U \cdot \n v_{-P}
+\int_{\frac{\O -Q}{\eps}} (V(\eps x +Q)-V(Q)) U v_{-P}
\\
=\int_{\frac{\O -Q}{\eps}} \left[-J(Q) \varDelta U + V(Q) U - U^p \right] v_{-P}
+J(Q) \int_{\de \OE} \frac{\de U_P}{\de \nu} v
\\
+\int_{\frac{\O -Q}{\eps}} (J(\eps x +Q)-J(Q)) \nabla U \cdot \n v_{-P}
+\int_{\frac{\O -Q}{\eps}} (V(\eps x +Q)-V(Q)) U v_{-P}.
\end{gather*}
Hence, since $U\equiv U^Q$ is solution of \eqref{eq:Q}, we get
\begin{multline}\label{eq:nf2}
(\nabla f_\eps(U_P) \mid v) =
J(Q) \int_{\de \OE} \frac{\de U_P}{\de \nu} v
+\int_{\frac{\O -Q}{\eps}} (J(\eps x +Q)-J(Q)) \nabla U \cdot \n v_{-P}
\\
+\int_{\frac{\O -Q}{\eps}} (V(\eps x +Q)-V(Q)) U v_{-P}.
\end{multline}
Let us estimate the first of these three terms:
\begin{gather*}
\left|J(Q) \int_{\de \OE} \frac{\de U_P}{\de \nu} v \right|
\le C \|v\|_{L^2(\de \OE)} 
\left( \int_{\de \OE}\left|\frac{\de U_P}{\de \nu}\right|^2\right)^{1/2}.
\end{gather*}
First of all, we observe that there exist $\eps_0>0$ and $C>0$ such that, for all 
$\eps \in(0,\eps_0)$ and for all $v\in H^1(\OE)$, we have
\[
\|v\|_{L^2(\de \OE)}\le C \|v\|_{H^1(\OE)}.
\]
Moreover, after making a translation and rotation, we can assume that $Q$ coincides with 
the origin $\cal O$ and that part of $\de \O$ is given by 
$x_N=\psi(x')=\frac 12 \sum^{N-1}_{i=1}\l_i x_i^2+O(|x'|^3)$ for $|x'| < \mu$, where 
$\mu$ is some constant depending only on $\O$. Then for $|y'|<\mu/\eps$, the corresponding 
part of $\de \OE$ is given by $y_N=\Psi (y')= \eps^{-1} \psi(\eps y')= 
\frac{\eps}{2} \sum^{N-1}_{i=1}\l_i y_i^2+O(\eps^2 |y'|^3)$. Then it is easy to see that
\[
\frac{\de U}{\de \nu}(y', \Psi(y'))
=\eps \left[\sum_{i=1}^{N-1} \l_i y_i \frac{\de U}{\de y_i}(y',0)
-\frac 12 \frac{\de^2 U}{\de y_N^2}(y',0) \sum_{i=1}^{N-1} \l_i y_i^2\right]
+O(\eps^2).
\]
Let us observe that by the exponential decay of $U$ and of its derivatives, we get:
\begin{gather*}
\int_{\de \tilde{\O}_\eps}\!\!\left|\frac{\de U}{\de \nu}\right|^2
\!\!\!= \! \eps^2 \!\int_{\de \tilde{\O}_\eps}\!\!
\left[\sum_{i=1}^{N-1} \l_i y_i \frac{\de U}{\de y_i}(y',0)
-\frac 12 \frac{\de^2 U}{\de y_N^2}(y',0) \sum_{i=1}^{N-1} \l_i y_i^2\right]^2
\!\!\!\!+\!o(\eps^2)
\!=\!O(\eps^2),
\end{gather*}
where $\de \tilde{\O}_\eps \equiv \de \OE \cap B_{\eps^{-1/2}}$. Therefore
\begin{equation}\label{eq:deU}
\left( \int_{\de \OE}\left|\frac{\de U}{\de \nu}\right|^2\right)^{1/2}
=\left( \int_{\de \OE \cap B_{\eps^{-1/2}}}\left|\frac{\de U}{\de \nu}\right|^2\right)^{1/2}
+o(\eps)=O(\eps).
\end{equation}
Let us calculate the second term of \eqref{eq:nf2}. We start observing that, from the 
assumption $D^2 J$ bounded, we infer that
\[
|J(\eps x +Q)-J(Q)| \le \eps |J'(Q)| |x| + c_1 \eps^2 |x|^2,
\]
and so, using again the exponential decay of $U$ and of its derivatives,
\begin{gather}
\int_{\frac{\O -Q}{\eps}} \!\!(J(\eps x +Q)-J(Q)) \nabla U \!\! \cdot \!\!\n v_{-P}
\le \|v\| \!\left(\! \int_{\frac{\O -Q}{\eps}}\!\! |J(\eps x +Q)-J(Q)|^2 |\nabla U|^2 \!\right)^{1/2} 
\nonumber
\\
\le c_2 \|v\| \left[\int_\RP \eps^2 |J'(Q)|^2 |x|^2|\nabla U|^4 
+ \int_\RP \eps^4 |x|^4 |\nabla U|^4 \right]^{1/2}
=O(\eps)\|v\|.  \label{eq:restoJ} 
\end{gather}
Analogously, we can say that:
\begin{equation}\label{eq:restoV}
\int_{\frac{\O -Q}{\eps}} (V(\eps x +Q)-V(Q)) U v_{-P}
=O(\eps)\|v\|.
\end{equation}
Now the conclusion follows immediately by \eqref{eq:nf2}, \eqref{eq:deU}, \eqref{eq:restoJ} 
and \eqref{eq:restoV}.
\end{proof}

We here present some useful estimates that will be used in the sequel.
 
\begin{proposition}\label{lemma1.2}
Let $P=Q/\eps \in \de \O_\eps$. Then we have:
\begin{equation}\label{eq:1.4}
\int_{\O_\eps} U_P^{p+1}=
\int_{\R^N_+}\left(U^Q\right)^{p+1} 
-\eps \frac{H(Q)}{2}
\int_{\R^{N-1}} \left[U^Q(y',0)\right]^{p+1} |y'|^2 d y'+o(\eps),
\end{equation}

\begin{equation}\label{eq:1.5}
\int_{\de \O_\eps} \frac{\de  U_P}{\de \nu} U_P
=-\eps\frac{(N-1)H(Q)}{4}  \int_{\R^{N-1}} \left[U^Q(y',0)\right]^2 d y'+o(\eps),
\end{equation}

\begin{multline}\label{eq:1.6}
J(Q) \int_{\O_\eps} |\n U_P|^2 
+V(Q) \int_{\O_\eps} U_P^2
\\
=\int_{\R^N_+}\left(U^Q\right)^{p+1} 
-\eps\frac{H(Q)}{2}\int_{\R^{N-1}} \left[U^Q(y',0)\right]^{p+1} |y'|^2 d y'
\\
-\eps J(Q) \frac{(N-1)H(Q)}{4} \int_{\R^{N-1}} \left[U^Q(y',0)\right]^2 d y'+o(\eps),     
\end{multline}

\begin{equation}\label{eq:J}
\int_{\O_\eps} J(\eps x)|\nabla U_P|^2=
J(Q) \int_{\OE}|\nabla U_P|^2
+\eps \int_{\R^-_{\nu(Q)}}J'(Q)[x] |\nabla U^Q|^2
+o(\eps),
\end{equation}

\begin{equation}\label{eq:V}
\int_{\O_\eps}  V(\eps x)  U_P^2=
V(Q) \int_{\OE}  U_P^2
+\eps \int_{\R^-_{\nu(Q)}} V'(Q)[x] \left(U^Q\right)^2 
+o(\eps).
\end{equation}
Moreover, we have
\begin{equation}\label{eq:1.8}
\int_{\O_\eps} U_P^p \; \de_{P_i} U_P=
\eps \frac{1}{p+1}  \bar C \de_{Q_i} \G(Q) + o(\eps),
\end{equation}

\begin{equation}\label{eq:1.7}
\de_{P_i} \left[ J(Q) \!\int_{\O_\eps} \!|\n U_P|^2 
+V(Q) \!\int_{\O_\eps}\! U_P^2 \right]
= \eps \bar C  \de_{Q_i} \G(Q) + o(\eps).
\end{equation}
where $\bar C= \int_\RP \bar U^{p+1}$ and $\G$ is defined in \eqref{eq:Gamma}.
\end{proposition}

\begin{proof}
The first two formulas can be proved repeating the arguments of 
Lemma 1.2 of \cite{Li}. Equation \eqref{eq:1.6} 
follows easily by \eqref{eq:1.4} and \eqref{eq:1.5} observing that
\[
J(Q) \int_{\O_\eps} |\n U_P|^2 
+V(Q) \int_{\O_\eps} U_P^2=
\int_{\OE}U_P^{p+1}
+J(Q)\int_{\de \OE}\frac{\de  U_P}{\de \nu} U_P.
\]
Let us prove \eqref{eq:J}. Arguing as in the proof of \eqref{eq:restoJ}, we infer:
\begin{gather*}
\int_{\O_\eps} J(\eps x)|\nabla U_P|^2=
\int_{\frac{\O -Q}{\eps}} J(\eps x +Q )|\nabla U^Q|^2
\\
=J(Q)\int_{\frac{\O -Q}{\eps}} |\nabla U^Q|^2
+\eps \int_{\frac{\O -Q}{\eps}} J'(Q)[x]|\nabla U^Q|^2 +o(\eps)
\\
=J(Q)\int_{\O_\eps} |\n U_P|^2
+\eps \int_{\R^-_{\nu(Q)}} J'(Q)[x]|\nabla U^Q|^2 
+o(\eps).
\end{gather*}
We can prove equation \eqref{eq:V} repeating the arguments of \eqref{eq:J}.

Since
\[ 
\int_{\O_\eps} U_P^p \; \de_{P_i} U_P=
\frac{1}{p+1} \de_{P_i} \int_{\O_\eps} U_P^{p+1},
\]
equations \eqref{eq:1.8} and \eqref{eq:1.7} follow easily because, as observed by \cite{Li}, the 
error terms $O(\eps )$ in \eqref{eq:1.4} and \eqref{eq:1.6} become of order $o(\eps)$ after applying 
$\de_{P_i}$ to them.
\end{proof}

\section{Invertibility of $D^2 f_\eps$ on $\left(T_{U_P}Z^\eps \right)^\perp$}

In this section we will show that $D^2 f_\eps$ is invertible on
$\left(T_{U_P}Z^\eps \right)^\perp$, where $T_{U_P} Z^\eps$ 
denotes the tangent space to $Z^\eps$ at $U_P$.

Let $L_{\eps,Q}:(T_{U_P}Z^\eps)^\perp\to
(T_{U_P}Z^\eps)^\perp$ denote the operator defined by setting
$(L_{\eps,Q}v \mid w)= D^2 f_\eps(U_P)[v,w]$.

\begin{lemma}\label{lem:inv}
There exists $C>0$ such that for $\eps$ small enough
one has that
\begin{equation}\label{eq:inv}
|(L_{\eps,Q}v \mid v)|\ge C \|v\|^{2},\qquad \forall\; v\in(T_{U_P}Z^{\eps})^{\perp}.
\end{equation}
\end{lemma}

\begin{proof}
By \eqref{eq:UQ}, if we set $\a(Q)=V(Q)^{\frac{1}{p-1}}$ and $\b(Q)=\sqrt{V(Q)/J(Q)}$, 
we have that $U^Q(x)=\a(Q) \bar U(\b(Q) x)$. Therefore, we have:
\begin{gather*}
\de_{P_i} U^Q(x-P)
= \de_{P_i} \left[\a(\eps P) \bar U(\b(\eps P)(x-P))\right] =
\\
\eps \de_{P_i} \a(\eps P) U^Q(\b(\eps P)(x-P))
\!+\! \eps \a(\eps P) \de_{P_i} \b(\eps P) \n U^Q(\b(\eps P)(x-P))\cdot (x-P)
\\
-\a(\eps P) \b(\eps P) (\de_{x_i} U^Q)(\b(\eps P)(x-P)).
\end{gather*}
Hence
\begin{equation} \label{eq:de_iU}
\de_{P_i} U^Q (x-P)=-\de_{x_i} U^Q(x-P)+O(\eps). 
\end{equation}

For simplicity, we can assume that $Q=\eps P$ is the origin $\cal O$.

Following \cite{Li}, without loss of generality, we assume that $Q=\eps P$ is 
the origin $\cal O$, $x_N$ is the tangent plane of 
$\de\O$ at $Q$ and $\nu(Q)=(0,\ldots, 0, -1)$. We also assume that part of $\de \O$ is given by 
$x_N=\psi(x')=\frac 12 \sum^{N-1}_{i=1}\l_i x_i^2+O(|x'|^3)$ for $|x'| < \mu$, where 
$\mu$ is some constant depending only on $\O$. Then for $|y'|<\mu/\eps$, the corresponding 
part of $\de \OE$ is given by $y_N=\Psi (y')= \eps^{-1} \psi(\eps y')= 
\frac{\eps}{2} \sum^{N-1}_{i=1}\l_i y_i^2+O(\eps^2 |y'|^3)$.

We recall that 
$T_{U^{\cal O}} Z^\eps = {\rm span}_{H^1(\OE)} \{\de_{P_1} U^{\cal O}, \ldots, \de_{P_{N-1}}U^{\cal O} \}$. 
We set
\begin{eqnarray*}
{\cal V}_\eps &=& 
{\rm span}_{H^1(\OE)} \{U^{\cal O}, \de_{x_1} U^{\cal O}, \ldots, \de_{x_{N-1}}U^{\cal O} \},
\\
{\cal V}_+ &=& 
{\rm span}_{H^1(\RP)} \{U^{\cal O}, \de_{x_1} U^{\cal O}, \ldots, \de_{x_{N-1}}U^{\cal O} \}.
\end{eqnarray*}
By \eqref{eq:de_iU} it suffices to prove (\ref{eq:inv}) for all 
$v\in {\rm span}\{U^{\cal O},\phi\}$, where $\phi$ is 
orthogonal to ${\cal V}_\eps$. Precisely we shall prove that there exist $C_{1},C_{2}>0$ such that,
for all $\eps>0$ small enough, one has:
\begin{eqnarray}
(L_{\eps,\cal O}U^{\cal O} \mid U^{\cal O})& \le & - C_{1}< 0.      \label{eq:neg} 
\\
(L_{\eps,\cal O}\phi \mid \phi)&\ge & C_{2} \|\phi\|^2.        \label{eq:claim}
\end{eqnarray}

The proof of (\ref{eq:neg}) follows easily from the fact that $U^{\cal O}$ is a Mountain Pass 
critical point of $F^{\cal O}$ and so from the fact that there exists $c_0>0$ 
such that, for all $\eps>0$ small enough, one finds:
\begin{equation*}
D^2 F^{\cal O}(U^{\cal O})[U^{\cal O},U^{\cal O}] < -c_0< 0.
\end{equation*}
Indeed, arguing as in the proof of Lemma \ref{eq:nf} (see \eqref{eq:restoJ} and \eqref{eq:restoV}) 
and by \eqref{eq:1.4} and \eqref{eq:1.6}, we have:
\begin{gather*}
(L_{\eps, \cal O} U^{\cal O} \mid U^{\cal O})=
\int_{\OE} J(\eps x) |\nabla U^{\cal O}|^2 
+ \int_{\OE} V(\eps x) (U^{\cal O})^2 
- p \int_{\OE} (U^{\cal O})^{p+1}
\\
=J({\cal O})\int_{\OE} |\nabla U^{\cal O}|^2 
+ V({\cal O}) \int_{\OE} (U^{\cal O})^2 
- p \int_{\OE} (U^{\cal O})^{p+1}+O(\eps)
\\
=D^2 F^{\cal O}(U^{\cal O})[U^{\cal O},U^{\cal O}]+O(\eps)< -c_0+O(\eps)<-C_1.
\end{gather*}

Let us prove (\ref{eq:claim}).

As before, the fact that $U^{\cal O}$ is a Mountain Pass critical point of $F^{\cal O}$ implies that
\begin{equation}\label{eq:D2F+}
D^2 F^{\cal O}(U^{\cal O})[\tilde \phi,\tilde \phi]>c_1 \|\tilde \phi\|^2_+ 
\quad \forall \tilde \phi \perp {\cal V}_+.
\end{equation}
Let us consider a smooth function
$\chi_{1}:\RN\to \R$ such that
\begin{equation*}
\chi_{1}(x) = 1, \quad \hbox{ for } |x| \le \eps^{-1/8}; \qquad
\chi_{1}(x) = 0, \quad \hbox{ for } |x| \ge 2 \eps^{-1/8};
\end{equation*}
\begin{equation*}
|\nabla \chi_{1}(x)| \le 2\eps^{1/8}, \quad \hbox{ for } \eps^{-1/8} \le |x| \le 2 \eps^{-1/8}.
\end{equation*}
We also set $ \chi_{2}(x)=1-\chi_{1}(x)$.
Given $\phi \perp {\cal V}_\eps$, let us consider the functions
\[
\phi_{i}(x)=\chi_{i}(x)\phi(x),\quad i=1,2.
\]
If $Q\neq {\cal O}$, then we would take
\[
\phi_{i}(x)=\chi_{i}(x-P)\phi(x),\quad i=1,2.
\]
With calculations similar to those of \cite{AMS}, we have
\begin{equation}\label{eq:phi}
\| \phi \|^2 = \| \phi_1 \|^2 + \| \phi_2 \|^2 + 
\underbrace{2\int_\RN \chi_{1}\chi_{2}(\phi^{2}+|\nabla \phi|^{2})}_{I_\phi} 
+ O(\eps^{1/8})\| \phi \|^2.
\end{equation}

We need to evaluate the three terms in the equation below:
\begin{equation}\label{eq:L}
(L_{\eps,{\cal O}}\phi \mid \phi)=
(L_{\eps,{\cal O}}\phi_{1} \mid \phi_{1})+
(L_{\eps,{\cal O}}\phi_{2} \mid \phi_{2})+
2(L_{\eps,{\cal O}}\phi_{1} \mid \phi_{2}).
\end{equation}
Let us start with $(L_{\eps,{\cal O}}\phi_{1} \mid \phi_{1})$.

Let $\eta=\eta_\eps$ a smooth cutoff function satisfying
\[
\eta(y) = 1, \quad \hbox{ for } |y| \le \eps^{-1/4}; \qquad
\eta(y) = 0, \quad \hbox{ for } |y| \ge 2 \eps^{-1/4};
\]
\[
|\nabla \eta(y)| \le 2\eps^{1/4}, \quad \hbox{ for } \eps^{-1/4} \le |y| \le 2 \eps^{-1/4}.
\]
Now we will straighten $\de \OE$ in the following way: let 
$\Phi \colon \RP \cap B_{\eps^{-1/2}} \to \OE$ be a function so defined:
\[
\Phi(y',y_N)=(y', y_N + \Psi(y')).
\]
We observe that:
\[
D \Phi (y)=
\left(
\begin{array}{ccc|c}
1   &                      &   & 0
\\
    &       \ddots         &   & \vdots 
\\
    &                      & 1 & 0
\\ \hline
    & \nabla_{y'} \Psi(y') &   & 1
\end{array}
\right).
\]
Let us defined $\tf\in H^1(\RP)$ as:
\[
\tf(y)=
\left\{
\begin{array}{lll}
\phi_1(\Phi(y))\, \eta(y) & \quad {\rm if }\;  |y| \le \eps^{-1/2},
\\
0 & \quad {\rm if }\; |y|>\eps^{-1/2}.
\end{array}
\right.
\]
We get:
\begin{gather*}
\int_{\RP}|\nabla \tf|^2 = 
\int_{\RP \cap B_{2 \eps^{-1/4}}}\left| \nabla \left[ \phi_1(\Phi(y)) \right]\right|^2 d y
\\
=\int_{\RP \cap B_{2 \eps^{-1/4}}} 
\sum_{i=1}^{N-1}\left| \frac{\de \phi_1}{\de x_i}(\Phi) 
+ \eps \l_i y_i \frac{\de \phi_1}{\de x_N}(\Phi)\right|^2
+\left|\frac{\de \phi_1}{\de x_N}(\Phi)\right|^2
+o(\eps)\|\phi\|^2
\\
=\int_{\RP \cap B_{2 \eps^{-1/4}}}|(\nabla  \phi_1)(\Phi)|^2  + O(\eps^{7/8})\|\phi\|^2= 
\int_{\OE}|\nabla  \phi_1|^2 + O(\eps^{7/8})\|\phi\|^2.
\end{gather*}
Analogously, we have:
\[
\int_{\RP} |\tf|^2 = \int_{\OE}|\phi_1|^2,
\]
and so
\[
\|\tf\|^2_+=\|\phi_1\|^2+ O(\eps^{7/8})\|\phi\|^2.
\]
Let us now evaluate $(L_{\eps,{\cal O}}\phi_{1}|\phi_{1})$:
\begin{gather*}
(L_{\eps,{\cal O}}\phi_{1} \mid \phi_{1})=
\int_{\OE} J(\eps x) |\nabla \phi_1|^2 
+\int_{\OE} V(\eps x) \phi_1^2 
-p \int_{\OE} (U^{\cal O})^{p-1} \phi_1^2
\\
=J({\cal O})\int_{\OE} |\nabla \phi_1|^2 
+V({\cal O})\int_{\OE} \phi_1^2 
-p \int_{\OE} (U^{\cal O})^{p-1} \phi_1^2
\\
+\eps \int_{\OE} J'({\cal O})[x] |\nabla \phi_1|^2 
+\eps \int_{\OE} V'({\cal O})[x] \phi_1^2 +o(\eps)\|\phi\|^2
\\
=J({\cal O})\int_{\OE} |\nabla \phi_1|^2 
+V({\cal O})\int_{\OE} \phi_1^2 
-p \int_{\OE} (U^{\cal O})^{p-1} \phi_1^2 +O(\eps^{7/8})\|\phi\|^2
\\
=J({\cal O})\int_{\RP} |\nabla \tf|^2 
+V({\cal O})\int_{\RP} \tf^2 
-p \int_{\RP} [U^{\cal O}(\Phi)]^{p-1} \tf^2 +O(\eps^{7/8})\|\phi\|^2
\\
=D^2 F^{\cal O}(U^{\cal O})[\tf,\tf]
-p \int_{\RP} \left( [U^{\cal O}(\Phi)]^{p-1}-(U^{\cal O})^{p-1}\right) \tf^2 +O(\eps^{7/8})\|\phi\|^2.
\end{gather*}
We have:
\begin{gather*}
\left|\int_{\RP} \left( [U^{\cal O}(\Phi)]^{p-1}-(U^{\cal O})^{p-1}\right) \tf^2 \right|
\le C \int_{\RP} |\Psi(y')| \tf^2 
\\
= O(\eps^{3/4})\|\tf\|^2 =O(\eps^{3/4})\|\phi\|^2.
\end{gather*}
Therefore, we have that
\begin{equation}\label{eq:L11}
(L_{\eps,{\cal O}}\phi_{1} \mid \phi_{1})=D^2 F^{\cal O}(U^{\cal O})[\tf,\tf] + O(\eps^{3/4})\|\phi\|^2.
\end{equation}
We can write $\tf=\xi + \zeta$, where $\xi \in {\cal V}_+$ and $\zeta \perp {\cal V}_+$. More precisely 
\[
\xi=(\tf  \mid  U^{\cal O})_+ \, U^{\cal O} \|U^{\cal O}\|^{-2}_+
+\sum_{i=1}^{N-1}(\tf  \mid  \de_{P_i} U^{\cal O})_+ \, \de_{P_i} U^{\cal O} \|\de_{P_i} U^{\cal O}\|^{-2}_+.
\]
Let us calculate $(\tf | U^{\cal O})_+$.
\begin{gather*}
(\tf  \mid  U^{\cal O})_+
=\int_{\RP} \nabla \tf \cdot \nabla  U^{\cal O}
+\int_{\RP} \tf  U^{\cal O}
\\
=\int_{\RP \cap B_{2 \eps^{-1/4}}}\nabla \left[ \phi_1(\Phi(y)) \right]\cdot \nabla  U^{\cal O}
+\int_{\RP \cap B_{2 \eps^{-1/4}}}  \phi_1(\Phi(y))\,  U^{\cal O}
\\
=\int_{\RP \cap B_{2 \eps^{-1/4}}} \!\!\!\!\!\!\!\left[(\nabla  \phi_1)(\Phi)\cdot \nabla  U^{\cal O} 
+ \phi_1(\Phi)\,  U^{\cal O}\right]
\!+\!\eps\! \sum_{i=1}^{N-1}\int_{\RP \cap B_{2 \eps^{-1/4}}} \!\!\!\!\!\!\!\!
\l_i y_i \frac{\de \phi_1}{\de x_N}(\Phi)  \frac{\de U^{\cal O}}{\de x_i}
\\
=\int_{\OE} \nabla \phi_1 \cdot \nabla  U^{\cal O}(\Phi^{-1})
+\int_{\OE} \phi_1  U^{\cal O}(\Phi^{-1})
+O(\eps^{7/8})\|\phi \|^2
\\
=\int_{\OE} \nabla \phi_1 \cdot \nabla  U^{\cal O}
+\int_{\OE} \phi_1  U^{\cal O}
+O(\eps^{3/4})\|\phi \|
=O(\eps^{3/4})\|\phi \|.
\end{gather*}
In an analogous way, we can prove also that $(\tf  \mid \de_{P_i} U^{\cal O})_+=O(\eps^{3/4})\|\phi \|$, 
and so
\begin{eqnarray}
\|\xi\|_+ &=&O(\eps^{3/4})\|\phi \|,    \label{eq:xi}
\\
\|\zeta\|_+ &=&\|\phi_1 \| + O(\eps^{3/4})\|\phi \|.     \label{eq:zeta}
\end{eqnarray}
Let us estimate $D^2 F^{\cal O}(U^{\cal O})[\tf,\tf]$. We get:
\begin{equation}\label{eq:D2F}
D^2 F^{\cal O}(U^{\cal O})[\tf,\tf]
=D^2 F^{\cal O}(U^{\cal O})[\zeta,\zeta]
+2 D^2 F^{\cal O}(U^{\cal O})[\zeta,\xi]
+D^2 F^{\cal O}(U^{\cal O})[\xi,\xi].
\end{equation}
By \eqref{eq:D2F+} and \eqref{eq:zeta}, we know that 
\[
D^2 F^{\cal O}(U^{\cal O})[\zeta,\zeta]>c_1 \|\zeta\|^2_+
=c_1 \|\phi_1 \|^2 + O(\eps^{3/4})\|\phi \|^2,
\]
while, by \eqref{eq:xi} and straightforward calculations, we have 
\begin{eqnarray*}
D^2 F^{\cal O}(U^{\cal O})[\zeta,\xi] &=& O(\eps^{3/4})\|\phi \|^2,
\\
D^2 F^{\cal O}(U^{\cal O})[\xi,\xi]  &=& O(\eps^{3/2})\|\phi \|^2.
\end{eqnarray*}
By these estimates, \eqref{eq:D2F} and \eqref{eq:L11}, we can say that
\begin{equation}\label{eq:L11-fin}
(L_{\eps,{\cal O}}\phi_{1} \mid \phi_{1})>c_1 \|\phi_1 \|^2 + O(\eps^{3/4})\|\phi\|^2.
\end{equation}

Using the definition of $\chi_i$ and the exponential decay of $U^{\cal O}$, we easily get
\begin{eqnarray}
(L_{\eps,{\cal O}}\phi_{2} \mid \phi_{2}) & \ge & c_2 \|\phi_{2}\|^{2}+o(\eps)\|\phi\|^{2},  \label{eq:L22}
\\
(L_{\eps,{\cal O}}\phi_{1} \mid \phi_{2}) & \ge & c_3 I_\phi + O(\eps^{1/8})\|\phi\|^{2},  \label{eq:L12}
\end{eqnarray}
where $I_\phi$ is defined in \eqref{eq:phi}.
Therefore by \eqref{eq:L}, \eqref{eq:L11-fin}, \eqref{eq:L22}, \eqref{eq:L12} and 
recalling \eqref{eq:phi} we get
\[
(L_{\eps,{\cal O}}\phi \mid \phi)\ge c_4 \|\phi\|^{2}+ O(\eps^{1/8})\|\phi\|^{2}.
\]
This completes the proof of the lemma.
\end{proof}

\section{The finite dimensional reduction}

\begin{lemma}\label{lem:w}
For $\eps>0$ small enough, there exists a unique
$w=w(\eps, Q)\in (T_{U_P} Z^\eps)^{\perp}$ such that
$\nabla f_\eps (U_P + w)\in T_{U_P} Z$.
Such a $w(\eps,Q)$ is of class $C^{2}$, resp.  $C^{1,p-1}$, with respect to $Q$, provided 
that $p\ge 2$, resp. $1<p<2$.
Moreover, the functional $\A_\eps (Q)=f_\eps (U_{Q/\eps} +w(\eps,Q))$ has
the same regularity of $w$ and satisfies:
\[
\nabla \A_\eps(Q_0)=0 
\quad \Longleftrightarrow \quad 
\nabla f_\eps\left(U_{Q_0/\eps}+w(\eps,Q_0)\right)=0.
\]
\end{lemma}

\begin{proof}
Let $\P=\P_{\eps, Q}$ denote the projection onto $(T_{U_P} Z^\eps)^\perp$. We want
to find a solution $w\in (T_{U_P} Z^\eps)^{\perp}$ of the equation
$\P\nabla f_\eps(U_P +w)=0$.  One has that $\nabla f_\eps(U_P+w)=
\nabla f_\eps (U_P)+D^2 f_\eps(U_P)[w]+R(U_P,w)$ with $\|R(U_P,w)\|=o(\|w\|)$, uniformly
with respect to $U_P$. Therefore, our equation is:
\begin{equation}\label{eq:eq-w}
L_{\eps,Q}w + \P\nabla f_\eps (U_P)+\P R(U_P,w)=0.
\end{equation}
According to Lemma \ref{lem:inv}, this is equivalent to
\[
w = N_{\eps,Q}(w), \quad \mbox{where}\quad
N_{\eps,Q}(w)=-L_{\eps,Q}\left( \P \nabla f_\eps (U_P)+\P R(U_P,w)\right).
\]
By \eqref{eq:nf} it follows that
\begin{equation}\label{eq:N}
\|N_{\eps,Q}(w)\| = O(\eps) + o(\|w\|).
\end{equation}
Then one readily checks that $N_{\eps,Q}$ is a contraction on some ball in
$(T_{U_P} Z^\eps)^{\perp}$
provided that $\eps>0$ is small enough.
Then there exists a unique $w$ such that $w=N_{\eps,Q}(w)$.  Let us
point out that we cannot use the Implicit Function Theorem to find
$w(\eps,Q)$, because the map $(\eps,u)\mapsto \P\nabla f_\eps (u)$ fails to be
$C^2$.  However, fixed $\eps>0$ small, we can apply the Implicit
Function Theorem to the map $(Q,w)\mapsto \P \nabla f_\eps (U_P + w)$.
Then, in particular, the function $w(\eps,Q)$ turns out to be of class
$C^1$ with respect to $Q$.  Finally, it is a standard argument, see
\cite{AB,ABC}, to check that the critical points of $\A_\eps(Q)=f_\eps (U_P+w)$ 
give rise to critical points of $f_\eps$.
\end{proof}

\begin{remark}\label{rem:w}
From (\ref{eq:N}) it immediately follows that:
\begin{equation}\label{eq:w}
\|w\|=O(\eps).
\end{equation}
\end{remark}

For future references, it is convenient to estimate the derivative $\de_{P_i} w$.

\begin{lemma}\label{lem:Dw}
If $\gamma=\min\{1,p-1\}$, then, for $i=1, \ldots, N-1$, one has that:
\begin{equation}\label{eq:Dw}
\|\de_{P_i} w\|=O(\eps^\gamma).
\end{equation}
\end{lemma}

\begin{proof}
We will set $h(U_P,w)=(U_P+w)^p -U_P^p -p U_P^{p-1}w$. With these notations, and recalling that
$L_{\eps,Q}w = -\dv (J(\eps x) \n w) +V(\eps x)w -p U_P^{p-1}w$, it follows that, 
for all $v\in (T_{U_P} Z^\eps)^{\perp}$, since $w$ satisfies \eqref{eq:eq-w}, then:
\begin{gather*}
\int_{\OE} J(\eps x) \n U_P \cdot \n v 
+\int_{\OE} V(\eps x) U_P v 
-\int_{\OE} U_P^p v
\\
+\!\int_{\OE} \!\!J(\eps x) \n w \cdot \n v 
+\!\int_{\OE} \!\!V(\eps x) w v 
-p\!\int_{\OE}\!\! U_P^{p-1} w v
- \!\int_{\OE}\!\! h(U_P,w)v =0.
\end{gather*}
Hence $\de_{P_i} w$ verifies:
\begin{gather}
\int_{\OE} J(\eps x) \n (\de_{P_i} U_P) \cdot \n v 
+\int_{\OE} V(\eps x) (\de_{P_i} U_P) v 
-p \int_{\OE} U_P^{p-1}(\de_{P_i} U_P)  v            \nonumber
\\
+\int_{\OE} J(\eps x) \n (\de_{P_i} w) \cdot \n v 
+\int_{\OE} V(\eps x) (\de_{P_i} w) v 
-p\int_{\OE} U_P^{p-1} (\de_{P_i} w) v            \nonumber
\\
-p(p-1)\int_{\OE} U_P^{p-2}(\de_{P_i} U_P)w v 
- \int_{\OE} \left[h_{U_P}(\de_{P_i} U_P) + h_w (\de_{P_i} w) \right]v =0.          \label{eq:dw}
\end{gather}
Let us set $L'=L_{\eps,Q}-h_w $. Then \eqref{eq:dw} can be written as
\begin{gather}
(L' (\de_{P_i} w)\mid v) 
=p(p-1)\int_{\OE} U_P^{p-2}(\de_{P_i} U_P)w v 
+ \int_{\OE} h_{U_P}(\de_{P_i} U_P) v        \nonumber
\\
-\int_{\OE} J(\eps x) \n (\de_{P_i} U_P) \cdot \n v 
-\int_{\OE} V(\eps x) (\de_{P_i} U_P) v 
+p \int_{\OE} U_P^{p-1}(\de_{P_i} U_P)  v.    \label{eq:L'}
\end{gather}
It is easy to see that 
\begin{equation}\label{eq:L'1}
\left| p(p-1)\int_{\OE} U_P^{p-2}(\de_{P_i} U_P)w v  \right| 
\le c_1 \|w\|\|v\| 
\end{equation}
and, if $\gamma = \min\{1,p-1\}$,
\begin{equation}\label{eq:L'2}
\left| \int_{\OE} h_{U_P}(\de_{P_i} U_P) v \right| 
\le c_2 \|w\|^\gamma \|v\|.
\end{equation}

Let us study the second line of \eqref{eq:L'}. We recall that often we will write $U$ instead of $U^Q$. 
Reasoning as in the proof of Lemma \ref{eq:nf} (see \eqref{eq:restoJ} and \eqref{eq:restoV}), 
we infer:
\begin{gather*}
I\equiv \int_{\OE} J(\eps x) \n (\de_{P_i} U_P) \cdot \n v 
+\int_{\OE} V(\eps x) (\de_{P_i} U_P) v 
-p \int_{\OE} U_P^{p-1}(\de_{P_i} U_P) v  
\\
=\int_{\frac{\O -Q}{\eps}} J(Q) \n (\de_{P_i} U) \cdot \n v_{-P}
+\int_{\frac{\O -Q}{\eps}} V(Q) (\de_{P_i} U) v_{-P}
\\
+\eps \int_{\OE} J'(Q)[x-P] \n (\de_{P_i} U_P) \cdot \n v
+\eps \int_{\OE} V'(Q)[x-P] (\de_{P_i} U_P) v
\\
-p \int_{\OE}  U_P^{p-1}(\de_{P_i} U_P)  v 
+O(\eps) \|v \|.
\end{gather*}

Suppose, for simplicity, $Q$ coincides with the origin $\cal O$ and that part of $\de \O$ is given by 
$x_N=\psi(x')=\frac 12 \sum^{N-1}_{i=1}\l_i x_i^2+O(|x'|^3)$ for $|x'| < \mu$, where 
$\mu$ is some constant depending only on $\O$. Then for $|y'|<\mu/\eps$, the corresponding 
part of $\de \OE$ is given by $y_N=\Psi (y')= \eps^{-1} \psi(\eps y')= 
\frac{\eps}{2} \sum^{N-1}_{i=1}\l_i y_i^2+O(\eps^2 |y'|^3)$. 

Since by \eqref{eq:de_iU} $\de_{P_i} U_P = -\de_{x_i} U_P+O(\eps)$, by integration by parts, we get:
\begin{eqnarray*}
\eps \int_{\OE} J'(Q)[x-P] \n (\de_{P_i} U_P)\! \cdot \!\n v
&=& \eps \int_{\OE} \de_{Q_i} J(Q) \n U_P \cdot \n v +O(\eps)\|v\|,
\\
\eps \int_{\OE} V'(Q)[x-P] (\de_{P_i} U_P) v
&=& \eps \int_{\OE} \de_{Q_i} V(Q) U_P v +O(\eps)\|v\|.
\end{eqnarray*}
Hence
\begin{gather*}
I =\int_{\OE} J(Q) \n (\de_{P_i} U_P) \cdot \n v
+\eps \int_{\OE} \de_{Q_i} J(Q) \n U_P \cdot \n v
\\
+\int_{\OE} V(Q) (\de_{P_i} U_P) v
+\eps \int_{\OE} \de_{Q_i} V(Q) U_P v
-p \int_{\OE}  U_P^{p-1}(\de_{P_i} U_P)  v 
+O(\eps) \|v \|.
\end{gather*}
Being $U=U^Q$ solution of \eqref{eq:Q}, we have that
\[
-J(Q) \varDelta (\de_{P_i} U) 
-\eps \de_{Q_i}J(Q) \varDelta U 
+V(Q) (\de_{P_i} U) 
+\eps \de_{Q_i} V(Q) U 
-p U^{p-1}(\de_{P_i} U)=0
\]
and so
\begin{gather*}
I=J(Q) \int_{\de \OE} \frac{\de}{\de \nu} (\de_{P_i} U_P)  v  
+\eps \de_{Q_i} J(Q) \int_{\de \OE} \frac{\de U_P}{\de \nu}   v  
+O(\eps) \|v \|. 
\end{gather*}
Arguing again as in the proof of Lemma \ref{eq:nf} (see \eqref{eq:deU}), 
we can prove that 
\[
\left|J(Q) \int_{\de \OE} \frac{\de}{\de \nu} (\de_{P_i} U_P)  v 
+\eps \de_{Q_i} J(Q) \int_{\de \OE} \frac{\de U_P}{\de \nu}   v\right|
= O(\eps) \|v \|. 
\]
Hence
\begin{equation}\label{eq:I}
I=O(\eps^{3/4}) \|v \|.
\end{equation}
Putting together \eqref{eq:L'}, \eqref{eq:L'1}, \eqref{eq:L'2} and \eqref{eq:I}, we find
\[
|(L'(\de w_i)\mid v)|=\left(c_3  \|w\|^{\gamma}+ O(\eps) \right)\|v \|.
\]
Since $h_w\to 0$ as $w \to 0$, the operator $L'$, likewise $L$, is 
invertible for $\eps >0$ small and therefore one finds
\[
\|\de_{P_i} w\|\le c_4  \|w\|^{\gamma}+ O(\eps).
\]
Finally, by Remark \ref{rem:w}, the Lemma follows.      
\end{proof}

\section{The finite dimensional functional}

\begin{theorem}\label{th:sviluppo}
Let $Q \in \de \O$ and $P=Q/\eps \in \de \OE$. Suppose {\bf (J)} and {\bf (V)}. Then, 
for $\eps$ sufficiently small, we get:
\begin{equation}\label{eq:A}
\A_\eps (Q)=
f_\eps(U_P + w(\eps,Q))
= c_0 \G(Q) +\eps \S(Q) +o(\eps),
\end{equation}
where $\G$ is the auxiliary functions introduced in \eqref{eq:Gamma},
\[
c_0\equiv\left(\frac{1}{2}-\frac{1}{p+1}\right)\int_{\RP} \bar U^{p+1},
\]
and $\S \colon \de \O \to \R$ is so defined:
\begin{multline}\label{eq:Sigma}
\S(Q)\equiv \frac{1}{2}\int_{\R^-_{\nu(Q)}} \!\!\!\!J'(Q)[x] |\nabla U^Q|^2 d x
+\frac{1}{2}\int_{\R^-_{\nu(Q)}}\!\!\!\!V'(Q)[x] \left(U^Q\right)^2 d x
\\
-\frac{1}{2} \bar{B}^Q J(Q) H(Q)
-\left(\frac{1}{2}-\frac{1}{p+1}\right) \bar{A}^Q H(Q),
\end{multline}
with
\begin{eqnarray*}
\bar{A}^Q &\equiv& \frac 12 \int_{\R^{N-1}}\left[U^Q(x',0)\right]^{p+1} |x'|^2 d x',  
\\
\bar{B}^Q &\equiv& \frac{(N-1)}{4} \int_{\R^{N-1}}\left[U^Q(x',0)\right]^2 d x.
\end{eqnarray*}
Moreover, for all $i=1,\ldots,N-1$, we get:
\begin{equation}\label{eq:DA}
\de_{P_i} \A_\eps (Q)= \eps c_0 \de_{Q_i} \G(Q)+o(\eps).
\end{equation}
\end{theorem}

\begin{proof}
In the sequel, to be short, we will often write $w$ instead of $w(\eps,Q)$. 
It is always understood that $\eps$ is taken in such a 
way that all the results discussed previously hold.

First of all, reasoning as in the proofs of \eqref{eq:J} and \eqref{eq:V} 
and by \eqref{eq:w}, we can observe that
\begin{eqnarray}
\int_{\OE} J(\eps x)\nabla U_P \cdot \nabla w 
&=&
J(Q) \int_{\OE} \nabla U_P \cdot \nabla w + o(\eps),     \label{eq:J2}
\\
\int_{\OE} V(\eps x)U_P \,w
&=& 
V(Q) \int_{\OE} U_P \,w +o(\eps).                        \label{eq:V2}
\end{eqnarray}
We have:
\begin{gather*}
\A_\eps (Q) =  
f_\eps(U_P + w(\eps,Q))
\\
=\frac{1}{2}\int_{\OE}\!\! J(\eps x) |\nabla (U_P +w)|^2
+\frac{1}{2}\int_{\OE}\!\! V(\eps x)(U_P+w)^2
-\frac{1}{p+1}\int_{\OE}\! (U_P+w)^{p+1}
\end{gather*}
[by \eqref{eq:w}]
\begin{gather*}
=\frac{1}{2}\int_{\OE} J(\eps x) |\nabla U_P|^2
+\frac{1}{2}\int_{\OE} V(\eps x)U_P^2
-\frac{1}{2}\int_{\OE} U_P^{p+1}
\\
+\int_{\OE} J(\eps x)\nabla U_P \cdot \nabla w 
+\int_{\OE} V(\eps x)U_P \,w
-\int_{\OE} U_P^p \,w
+\left(\frac{1}{2}-\frac{1}{p+1}\right) \int_{\OE} \!\!U_P^{p+1}
\\
-\frac{1}{p+1}\int_{\OE}\!\left[ (U_P+w)^{p+1} -U_P^{p+1}-(p+1) U_P^p \,w\right]
+o(\eps)=
\end{gather*}
[by \eqref{eq:1.6}, \eqref{eq:J}, \eqref{eq:V}, \eqref{eq:J2} and \eqref{eq:V2}
and with our notations]
\begin{gather*}
=\frac{1}{2} \int_{\RP} U^{p+1}
-\frac{\eps}{2} \bar{A}^Q  H(Q)
-\frac{\eps}{2} \bar{B}^Q  J(Q) H(Q)
+\frac{\eps}{2}\int_{\R^-_{\nu(Q)}}\!\!\!\!J'(Q)[x] |\nabla U|^2
\\
+\frac{\eps}{2}\int_{\R^-_{\nu(Q)}}V'(Q)[x] U^2
-\frac{1}{2}\int_{\RP} U^{p+1}
+\frac{\eps}{2} \bar{A}^Q H(Q)
\\
+J(Q) \int_{\OE}  \nabla U_P \cdot \nabla w 
+V(Q) \int_{\OE}  U_P \,w
-\int_{\OE} U_P^p \,w
\\
+\left(\frac{1}{2}-\frac{1}{p+1}\right) \int_{\RP} U^{p+1}
-\eps \left(\frac{1}{2}-\frac{1}{p+1}\right) \bar{A}^Q H(Q)
+o(\eps).
\end{gather*}
From the fact that $U$ is solution of \eqref{eq:Q}, we infer
\begin{gather*}
J(Q) \int_{\OE}  \nabla U_P \cdot \nabla w 
+V(Q) \int_{\OE}  U_P \,w
-\int_{\OE} U_P^p \,w
\\
=\int_{\OE} \left[ -J(Q) \varDelta U_P + V(Q)U_P - U_P^p \right] w 
+J(Q)\int_{\de \OE} \frac{\de U_P}{\de \nu} w
\\
=J(Q)\int_{\de \OE} \frac{\de U_P}{\de \nu} w
=o(\eps).
\end{gather*}
By these considerations we can say that
\begin{multline*}
\A_\eps (Q) = 
\left(\frac{1}{2}-\frac{1}{p+1}\right) \int_{\RP} U^{p+1}
\\
+\eps \Bigg[
\frac{1}{2}\int_{\R^-_{\nu(Q)}} J'(Q)[x] |\nabla U|^2
+\frac{1}{2}\int_{\R^-_{\nu(Q)}}V'(Q)[x] U^2
\\
-\frac{1}{2} \bar{B}^Q J(Q) H(Q)
-\left(\frac{1}{2}-\frac{1}{p+1}\right) \bar{A}^Q H(Q)
\Bigg]
+o(\eps).
\end{multline*}
Now the conclusion of the first part of the theorem follows observing that, since by \eqref{eq:UQ}
\[
U^Q (x)=V(Q)^{\frac{1}{p-1}}\,\bar U\left(x \sqrt{V(Q)/J(Q)}  \right),
\]
then
\[
\int_{\RP} U^{p+1}= V(Q)^{\frac{p+1}{p-1}-\frac{N}{2}}J(Q)^{\frac{N}{2}} \int_{\RP} \bar U^{p+1}.
\]

Let us prove now the estimate on the derivatives of $\A_\eps$. First of all, we observe that 
by \eqref{eq:nf} and by \eqref{eq:Dw}, we infer that
\[
\left|\n f_\eps (U_P)[\de_{P_i} w]\right| =O(\eps^{1+\g}),
\]
and so, by \eqref{eq:w} and \eqref{eq:Dw}, we have:
\begin{gather*}
\de_{P_i} \A_\eps (Q)
=\n f_\eps (U_P +w)[\de_{P_i} U_P +\de_{P_i} w]=
\n f_\eps (U_P +w)[\de_{P_i} U_P]+ O(\eps^{1+\g})
\\
=\n f_\eps (U_P)[\de_{P_i} U_P]
+D^2  f_\eps (U_P) [w, \de_{P_i} U_P]
\\
+\left(\n f_\eps (U_P +w)
-\n f_\eps (U_P)
-D^2 f_\eps (U_P)[w] \right) [\de_{P_i} U_P] + O(\eps^{1+\g}).
\end{gather*}
But 
\[
\| \n f_\eps (U_P +w)
-\n f_\eps (U_P)
-D^2 f_\eps (U_P)[w] \|
= o(\|w \|) =o(\eps)
\]
and, moreover, by \eqref{eq:eq-w} also $D^2  f_\eps (U_P) [w, \de_{P_i} U_P] =O(\eps^{1+\g})$, therefore
\begin{equation}\label{eq:DA-1}
\de_{P_i} \A_\eps (Q)
=\n f_\eps (U_P)[\de_{P_i} U_P] + O(\eps^{1+\g}).
\end{equation}
Let us calculate $\n f_\eps (U_P)[\de_{P_i} U_P]$.
\begin{gather*}
\n f_\eps (U_P)[\de_{P_i} U_P]
\!=\!\!\int_{\OE}\!\!\!\! J(\eps x) \n U_P\! \cdot \!\n (\de_{P_i} U_P)
\!+\!\!\int_{\OE}\!\!\!\! V(\eps x) U_P (\de_{P_i} U_P)
\!-\!\!\int_{\OE}\!\!\!\! U_P^p  (\de_{P_i} U_P)
\\
=J(Q) \!\int_{\OE}\!\! \n U_P \cdot \n (\de_{P_i} U_P)
+V(Q) \!\int_{\OE } \!\!U_P (\de_{P_i} U_P)
\\
+\eps \!\!\int_{\R^-_{\nu(Q)}}\!\!\!\!\!\!\!\!J'(Q)[x] \n U\! \cdot \!\n (\de_{P_i} U)\!
+\!\eps\!\! \int_{\R^-_{\nu(Q)}}\!\!\!\!\!\!\!\!V'(Q)[x] U (\de_{P_i} U)\!
-\!\!\int_{\OE}\!\!\!\!  U_P^p  (\de_{P_i} U_P)
\!+\!o(\eps).
\end{gather*}

Suppose, for simplicity, $Q$ coincides with the origin $\cal O$ and that part of $\de \O$ is given by 
$x_N=\psi(x')=\frac 12 \sum^{N-1}_{i=1}\l_i x_i^2+O(|x'|^3)$ for $|x'| < \mu$, where 
$\mu$ is some constant depending only on $\O$. Then for $|y'|<\mu/\eps$, the corresponding 
part of $\de \OE$ is given by $y_N=\Psi (y')= \eps^{-1} \psi(\eps y')= 
\frac{\eps}{2} \sum^{N-1}_{i=1}\l_i y_i^2+O(\eps^2 |y'|^3)$. 

Since by \eqref{eq:de_iU} $\de_{P_i} U_P = -\de_{x_i} U_P +O(\eps)$, by integration by parts, we get:
\begin{eqnarray*}
\int_{\R^-_{\nu(Q)}}\!\!\!J'(Q)[x] \n U\! \cdot \!\n (\de_{P_i} U)
&=&\frac 12 \int_{\R^-_{\nu(Q)}}\! \de_{Q_i} J(Q) |\n U|^2,
\\
\int_{\R^-_{\nu(Q)}}\!V'(Q)[x] U (\de_{P_i} U)
&=&\frac 12 \int_{\R^-_{\nu(Q)}}\! \de_{Q_i} V(Q) U^2.
\end{eqnarray*}
Therefore we infer 
\[
\n f_\eps (U_P)[\de_{P_i} U_P]
=\frac 12 \de_{P_i}\!\!\left[ J(Q) \!\!\int_{\OE}\!\!|\n U_P|^2\!
+\!V(Q) \!\!\int_{\OE} \!\!U_P^2 \right]
-\int_{\OE} \!\!U_P^p \; (\de_{P_i} U_P)+o(\eps),
\]
and so, by \eqref{eq:1.8} and \eqref{eq:1.7},
\[
\n f_\eps (U_P)[\de_{P_i} U_P]
= \eps \left[\left( \frac 12 - \frac{1}{p+1}\right) \int_\RP \bar U^{p+1}\right] \de_{Q_i} \G(Q) 
= \eps c_0 \de_{Q_i} \G(Q)+o(\eps).
\]
By this equation and by \eqref{eq:DA-1}, \eqref{eq:DA} follows immediately.
\end{proof}

\begin{remark}\label{re:Gamma}
Let us observe that by \eqref{eq:A} and \eqref{eq:DA}, for $\eps$ sufficiently small,
we have
\begin{equation}\label{eq:C1}
\|\A_\eps - c_0 \G \|_{C^1(\de \O)}=O(\eps).
\end{equation}
\end{remark}

\begin{remark}\label{re:Sigma}
By \eqref{eq:UQ}, it is easy to see that, if $J$ and $V$ are constant on the boundary $\de \O$, 
then $\bar \S$, defined in \eqref{eq:Sigmabar}, coincides with $\S$, defined in \eqref{eq:Sigma} 
with the following definitions:
\begin{equation*}
\begin{array}{lll}
C_J \equiv J_{|_{\de \O}},& \quad&C_V \equiv V_{|_{\de \O}}, 
\\
k_1 \equiv \frac{(C_V)^{\frac{p+1}{p-1}}}{2 C_J},& &k_2 \equiv \sqrt{C_V/C_J},
\\
k_3 \equiv \frac{(C_V)^{\frac{2}{p-1}}}{2}, & &k_4 \equiv -\frac{1}{2} \bar{B} C_J 
-\left(\frac{1}{2}-\frac{1}{p+1}\right) \bar{A}, 
\end{array}
\end{equation*}
where 
\begin{eqnarray*}
\bar{A} &\equiv& 
\frac {(C_V)^{\frac{p+1}{p-1}}}{2} 
\int_{\R^{N-1}}\left[\bar U \left(x'\sqrt{C_V/C_J},0 \right)\right]^{p+1} |x'|^2 d x',
\\
\bar{B} &\equiv& 
\frac{(N-1)(C_V)^{\frac{2}{p-1}}}{4} 
\int_{\R^{N-1}}\left[\bar U \left(x'\sqrt{C_V/C_J},0 \right)\right]^2 d x'.   
\end{eqnarray*}
\end{remark}

\section{Proofs of Theorem \ref{th1} and Theorem \ref{th2}}

In this section we will state and prove two multiplicity results for \eqref{EQe} whose 
Theorem \ref{th1} is a particular case. Finally we will prove also 
Theorem \ref{th2}.

Let us start introducing a topological invariant related to Conley theory.

\begin{definition}
Let $M$ be a subset of $\RN$, $M\ne \emptyset$.
The {\sl cup long} $l(M)$ of $M$ is defined by
$$
l(M)=1+\sup\{k\in \N \mid \exists\, \a_{1},\ldots,\a_{k}\in
\check{H}^{*}(M)\setminus 1, \,\a_{1}\cup\ldots\cup\a_{k}\ne 0\}.
$$
If no such class exists, we set $l(M)=1$. Here $\check{H}^{*}(M)$ is the
Alexander
cohomology of $M$ with real coefficients and $\cup$ denotes the cup product.
\end{definition}

Let us recall Theorem 6.4 in Chapter II of \cite{C}.

\begin{theorem}\label{th:Chang}
Let $N$ a Hilbert-Riemannian manifold. 
Let $g \in C^2(N)$ and let $M\subset N$ be a smooth compact nondegenerate manifold of 
critical points of $g$. Let $U$ be a neighborhood of $M$ and let $h \in C^1(N)$. 
Then, if $\|g-h\|_{C^1(\bar{U})}$ is sufficiently small, the function $g$ possesses 
at least $l(M)$ critical points in $U$.
\end{theorem}

Let us suppose that $\G$ has a smooth manifold of critical points $M$. We say that $M$
is nondegenerate (for $\G$) if every $x\in M$ is a nondegenerate
critical point of $\G_{|M^{\perp}}$. The Morse index of $M$ is, by definition,
the Morse index of any $x\in M$, as critical point of $\G_{|M^{\perp}}$.

We now can state our first multiplicity result.

\begin{theorem}\label{th:lc}
Let {\bf (J)} and {\bf (V)} hold and suppose $\G$ has a nondegenerate smooth manifold of critical points 
$M\subset \de \O$. 
There exists $\eps_0>0$ such that if $0<\eps<\eps_0$, then \eqref{EQe} has at least $l(M)$ solutions that
concentrate near points of $M$.
\end{theorem}

\begin{proof}
Fix a $\delta$-neighborhood $M_\delta$ of $M$ such that 
the only critical points of $\G$ in $M_\delta$ are those in $M$. We will take $U=M_\delta$. 

For $\eps$ sufficiently small, by \eqref{eq:C1} and Theorem \ref{th:Chang}, 
$\A_\eps$ possesses at least $l(M)$ critical points, which are solutions of \eqref{EQ} 
by Lemma \ref{lem:w}. Let $Q_\eps \in M$ be one of these critical points, then 
$u_\eps^{Q_\eps}=U_{Q_\eps/\eps}+w(\eps, Q_\eps)$ is a solution of \eqref{EQ}. Therefore 
\[
u_\eps^{Q_\eps}(x/\eps)\simeq U_{Q_\eps/\eps}(x/\eps)= U^{Q_\eps}\left(\frac{x-Q_\eps}{\eps} \right)
\]
is a solution of \eqref{EQe}.
\end{proof}

Moreover, when we deal with local minima (resp. maxima) of $\G$, the
preceding results can be improved because the number of positive solutions of \eqref{EQe} 
can be estimated by means of the category and $M$ does not need to be a manifold.

\begin{theorem}\label{th:cat}
Let {\bf (J)} and {\bf (V)} hold and suppose $\G$ has
a compact set $X \subset \de \O$ where $\G$ achieves a strict local minimum (resp. maximum), 
in the sense that there exist $\delta>0$ and a $\d$-neighborhood $X_\d \subset \de \O$ of $X$ such that
\[
b\equiv \inf\{\G(Q) : Q\in \partial X_{\d}\}> a\equiv \G_{|_X}, \quad
\left({\rm resp. }\; \sup\{\G(Q) : Q\in \partial X_{\d}\}< \G_{|_X} \right).
\]

Then there exists $\eps_0>0$ such that \eqref{EQe}
has at least $\cat(X,X_\d)$ solutions that concentrate near points of $X_{\d}$, provided 
$\eps \in (0,\eps_0)$. Here $\cat(X,X_\d)$ denotes the Lusternik-Schnirelman category of $X$ with 
respect to $X_\d$. 
\end{theorem}

\begin{proof}
We will treat only the case of minima, being the other one similar.
We set 
%
%
$Y=\{Q\in X_{\d} :\A_{\eps}(Q)\le c_0 (a+b)/2\}$.
By \eqref{eq:A} it follows that there exists $\eps_0>0$ such that
\begin{equation}\label{eq:X}
X  \subset Y \subset X_{\d},
\end{equation}
provided $\eps\in (0,\eps_0)$. Moreover, if $Q \in \partial X_{\d}$ then
$\G(Q)\ge b$ and hence
\[
\A_{\eps}(Q)\ge c_0 \G (Q) + O(\eps) \ge  c_0 b + O(\eps).
\]
On the other side, if $Q\in Y$
then $\A_{\eps}(Q)\le c_0 (a+ b)/2$.  Hence, for $\eps$ small, 
$Y$ cannot meet $\partial X_{\d}$ 
and this readily implies that $Y$ is compact.
Then $\A_{\eps}$ possesses at least $\cat(Y,X_{\d})$
critical points in $ X_{\d}$.  Using (\ref{eq:X}) and the properties of
the category one gets
\[
\cat(Y,Y)\ge \cat(X,X_{\d}),
\]
and the result follows.
\end{proof}

\begin{remark}
Let us observe that the (a) of Theorem \ref{th1} is a particular case of Theorem \ref{th:lc} while  
the (b) of Theorem \ref{th1} is a particular case of Theorem \ref{th:cat}.
\end{remark}

\noindent Let us now prove Theorem \ref{th2}.

\begin{proof2}
Let $Q$ be a minimum point of $\bar \S$ (the other case is similar) and let 
$\L \subset \de \O$ be a compact neighborhood of $Q$ such that 
\begin{equation*}
\min_\L \bar \S<\min_{\de \L}\bar \S.
\end{equation*}
By \eqref{eq:A} and Remark \ref{re:Sigma}, it is easy to see that for $\eps$ 
sufficiently small, there results:
\[
\min_{\L} \A_\eps <\min_{\de \L}\A_\eps.
\]
Hence, $\A_\eps$ possesses a critical point $Q_\eps$ in $\L$. By Lemma \ref{lem:w} we have that 
$u_{\eps, Q_\eps}=U_{Q_\eps/\eps}+ w(\eps, Q_\eps)$ is a critical point of $f_\eps$ and so a solution of 
problem \eqref{EQ}. Therefore
\[
u_{\eps, Q_\eps}(x/\eps)\simeq U_{Q_\eps/\eps}(x/\eps)= U^{Q_\eps}\left(\frac{x-Q_\eps}{\eps} \right)
\]
is a solution of (\ref{EQe}).
\end{proof2}


\begin{thebibliography}{99}

\bibitem{AB}
A. Ambrosetti, M. Badiale,
{\it Variational perturbative methods and bifurcation of bound states from the essential spectrum}, 
Proc. Royal Soc. Edinburgh, {\bf 128 A}, (1998), 1131--1161.

\bibitem{ABC}
A. Ambrosetti, M. Badiale, S. Cingolani, 
{\it Semiclassical states of nonlinear Schr\"odinger equations}, 
Arch. Rat. Mech. Anal., {\bf 159}, (2001), 253--271.


\bibitem{AMS}
A. Ambrosetti, A. Malchiodi, S. Secchi, 
{\it Multiplicity results  for some nonlinear Schr\"odinger equations with potentials}, 
Arch. Rational Mech. Anal., {\bf 159}, (2001), 253--271.

\bibitem{C}
K. C. Chang, Infinite dimensional Morse theory and multiple solutions problems, Birkh\"auser, 1993.

\bibitem{DY}
E. N. Dancer, S. Yan, {\it Multipeak solutions for a singularly perturbed Neumann problem}, 
Pacific J. Math., {\bf 189}, no. 2, (1999), 241--262.



\bibitem{dPFW}
M. del Pino, P. Felmer, J. Wei, 
{\it On the role of mean curvature in some singularly perturbed Neumann problems}, 
SIAM J. Math. Anal., {\bf 31}, (1999), no. 1, 63--79.

\bibitem{GPW}
M. Grossi, A. Pistoia, J. Wei, {\it Existence of multipeak solutions for a semilinear Neumann problem via 
nonsmooth critical point theory}, 
Cal. Var. PDE, {\bf 11}, no. 2, (2000), 143--175.

\bibitem{G}
C. Gui, {\it Multipeak solutions for a semilinear Neumann problem}, 
Duke Math. J., {\bf 84}, no. 3, (1996), 739--769.

\bibitem{GW}
C. Gui, J. Wei, 
{\it On multiple mixed interior and boundary peak solutions for some singularly perturbed Neumann problems}, 
Canad. J. Math., {\bf 52}, no. 3, (2000), 522--538.

\bibitem{Li}
Y. Y. Li, {\it On a singularly perturbed equation with Neumann boundary condition}, 
Comm. PDE, {\bf 23(3\&4)}, (1998), 487--545.  

\bibitem{NT}
W. M. Ni, I. Takagi, 
{\it On the shape of the least-energy solutions to a semilinear Neumann problem}, 
Comm. Pure Appl. Math., {\bf 44}, (1991), 819--851.

\bibitem{NT2}
W. M. Ni, I. Takagi, 
{\it Locating the peaks of least-energy solutions to a semilinear Neumann problem}, 
Duke Math. J., {\bf 70}, (1993), 247--281.


\bibitem{PS}
A. Pomponio, S. Secchi, 
{\it On a class of singularly perturbed elliptic equations in divergence form: 
existence and multiplicity results}, Preprint SISSA, 36/2003/M, (2003).


\bibitem{W2}
J. Wei, {\it On the boundary spike layer solutions to a singularly perturbed Neumann problem}, 
J. Diff. Eq., {\bf 134}, no. 1, (1997), 104--133.

\end{thebibliography}
\end{document}